\documentclass{amsart}

\newtheorem{theorem}{Theorem}[section]
\newtheorem{lemma}[theorem]{Lemma}

\theoremstyle{definition}
\newtheorem{definition}[theorem]{Definition}

\newtheorem{proposition}[theorem]{Proposition}
\newtheorem{corollary}[theorem]{Corollary}
\theoremstyle{remark}
\newtheorem{remark}[theorem]{Remark}

\numberwithin{equation}{section}

\begin{document}

\title[On a two-variable $p$-adic $l_q$-function]
{On a two-variable $p$-adic $l_q$-function}

\author{Min-Soo Kim, Taekyun Kim, D. K. Park and Jin-Woo Son}

\address{Department of Mathematics \\ Kyungnam University \\ Masan 631-701, Korea}
\email{mskim@kyungnam.ac.kr}

\address{EECS, Kyungpook National University \\ Taegu 702-701, Korea}
\email{tkim@knu.ac.kr}

\address{Department of Physics \\ Kyungnam University \\ Masan 631-701, Korea}
\email{dkpark@kyungnam.ac.kr}

\address{Department of Mathematics \\ Kyungnam University \\ Masan 631-701, Korea}
\email{sonjin@kyungnam.ac.kr}



\begin{abstract}
We prove that a two-variable $p$-adic $l_q$-function has the series expansion
$$l_{p,q}(s,t,\chi)=\frac{[2]_q}{[2]_{q^F}}
\sum_{\substack{a=1\\ (p,a)=1}}^F(-1)^a\frac{\chi(a)q^a}{\langle a+pt\rangle^{s}}\sum_{m=0}^\infty\binom {-s}m
\left(\frac F{\langle a+pt\rangle}\right)^{m}
E^*_{m,q^F}$$
which interpolates a linear combinations of terms of the generalized $q$-Euler polynomials at non positive integers.
The proof of this original construction is due to Kubota and Leopoldt in 1964, although the method given this note is due to Washington.
\end{abstract}

\maketitle

\section{Introduction}

The ordinary Euler polynomials $E_{n}(t)$ are defined by the equation
$$\frac{2e^{tx}}{e^x+1}=\sum_{n=0}^\infty E_n(t)\frac{x^n}{n!}.$$
Setting $t=1/2$ and normalizing by $2^n$ gives the ordinary Euler numbers
$$E_n=2^nE_n\left(\frac12\right).$$
The ordinary Euler polynomials appear in many classical results (see \cite{AS}).
In \cite{CK}, the values of these polynomials at rational arguments were expressed in term of the Hurwitz zeta function.
Congruences for Euler numbers have also received much attention from the point of view of $p$-adic interpolation.
In \cite{KKJR}, Kim {\it et al.} recently defined the natural $q$-extension of ordinary Euler numbers and polynomials
by $p$-adic integral representation and proved properties generalizing those satisfied
by $E_n$ and $E_n(t).$  They also constructed the one-variable $p$-adic $q$-$l$-function $l_{p,q}(s,\chi)$
for Dirichlet characters $\chi$ and $s\in\mathbb C_p$ with $|s|_p<p^{1-\frac1{p-1}},$ with the property that
$$l_{p,q}(-n,\chi)=E^*_{n,\chi\omega^{-n},q}-{[2]_q}{[2]_{q^p}^{-1}}p^n\chi\omega^{-n}(p)E^*_{n,\chi\omega^{-n},q^p}$$
for $n=0,1,\ldots,$ where $E^*_{n,\chi\omega^{-n},q}$ is a generalized $q$-Euler numbers attached to the
Dirichlet characters $\chi\omega^{-n}$ (see Section \ref{def} for definitions).

In the present paper, we shall construct a specific two-variable $p$-adic $l_q$-function $l_{p,q}(s,t,\chi)$
by means of a method provided in \cite{Wa2, Fo, KT4}.
We also prove that $l_{p,q}(s,t,\chi)$ is analytic in $s$ and $t$ for $s\in\mathbb C_p$ with $|s|_p<p^{1-\frac1{p-1}}$ and
$t\in\mathbb C_p$ with $|t|_p\leq1,$ which interpolates a linear combinations of terms of the generalized $q$-Euler polynomials at non positive integers.
This two-variable function is a generalization
of the one-variable $p$-adic $q$-$l$-function, which is the function obtained by putting $t=0$ in $l_{p,q}(s,t,\chi)$
(cf. \cite{Fo, KT1, KT2, KT3, KT4, KKJR, Yo, Wa1, Wa2}).

Thought this paper $\mathbb Z, \mathbb Z_p, \mathbb Q_p$ and $\mathbb C_p$
will denote the ring of integers, the ring of $p$-adic rational
integers, the field of $p$-adic rational numbers and completion
of the algebraic closure of $\mathbb Q_p,$ respectively.
We will use $\mathbb Z^+$ for the set of non positive integers.
Let $v_p$ be the normalized exponential valuation of $\mathbb C_p$ with
$|p|_p=\frac1p.$ When one talks of $q$-extension, $q$ is
variously considered as an indeterminate, a complex number $q\in\mathbb
C,$ or a $p$-adic number $q\in\mathbb C_p.$ If $q\in\mathbb C_p,$ then we
normally assume $|1-q|_p<1.$ If $q\in\mathbb C,$ then we assume that
$|q|<1.$ Also we use the following notations:
$$[x]_q=\frac{1-q^x}{1-q}\quad\text{and}\quad [x]_{-q}=\frac{1-(-q)^x}{1+q},\quad \text{cf. \cite{KT1, KT2}}.$$

Let $d$ be a fixed integer, and let
$$\begin{aligned}
&X=X_d=\varprojlim_N (\mathbb Z/dp^N\mathbb Z),\quad
X^*=\bigcup_{\substack{ 0<a<dp\\ (a,p)=1}} a+dp\mathbb Z_p,\\
&a+dp^N\mathbb Z_p=\{x\in X\mid x\equiv a\pmod{dp^N}\},
\end{aligned}$$
where $a\in \mathbb Z$ lies in $0\leq a<dp^N$. Let $UD(\mathbb Z_p)$ be the space of uniformly differentiable function on
$\mathbb Z_p.$ For $f\in UD(\mathbb Z_p),$ the $p$-adic $q$-integral was defined by
\begin{equation}\label{I_q}
\begin{aligned}
I_q(f)&=\int_{\mathbb Z_p}f(a)d\mu_q(a)=\int_Xf(a)d\mu_q(a) \\
&=\lim_{N\rightarrow\infty}\frac1{[dp^N]_q}\sum_{a=0}^{dp^N-1}f(a)q^a\quad\text{for } |1-q|_p<1.
\end{aligned}
\end{equation}
In \cite{KT2}, the bosonic integral was considered from a more
physical point of view to the bosonic limit $q\rightarrow1$ as
follows:
\begin{equation}\label{I_1}
I_1(f)=\lim_{q\rightarrow1} I_q(f)=\int_{\mathbb Z_p}f(a)d\mu_1(a)=\lim_{N\rightarrow\infty}\frac1{p^N}
\sum_{a=0}^{p^N-1}f(a).
\end{equation}
Furthermore, we can consider the fermionic integral in contrast to the conventional ``bosonic.''
That is,
$I_{-1}(f)=\int_{\mathbb Z_p}f(a)d\mu_{-1}(a)$ (see \cite{KT3}).
From this, we derive
$I_{-1}(f_1)+I_{-1}(f)=2f(0),$
where $f_1(a)=f(a+1).$ Also we have
\begin{equation}\label{I_-1}
I_{-1}(f_n)+(-1)^{n-1}I_{-1}(f)=2\sum_{a=0}^{n-1}(-1)^{n-1-a}f(a),
\end{equation}
where $f_n(a)=f(a+n)$ and $n\in\mathbb Z^+$ (see \cite{KT3}).
For $|1-q|_p<1,$ we consider fermionic $p$-adic $q$-integral on $\mathbb Z_p$ which is the $q$-extension of $I_{-1}(f)$
as follows:
\begin{equation}\label{I_-q}
I_{-q}(f)=\int_{\mathbb Z_p}f(a)d\mu_{-q}(a)
=\lim_{N\rightarrow\infty}\frac1{[dp^N]_{-q}}\sum_{a=0}^{dp^N-1}f(a)(-q)^a\quad
\text{(cf. \cite{KKJR}). }
\end{equation}

\section{$q$-Euler numbers and polynomials}\label{def}

In this section, we review some notations and facts in \cite{KKJR}.

From (\ref{I_-q}), we can derive the following formula:
\begin{equation}\label{f-e}
qI_{-q}(f_1)+I_{-q}(f)=[2]_qf(0),
\end{equation}
where $f_1(a)$ is translation with $f_1(a)=f(a+1).$ If we take
$f(a)=e^{ax},$ then we have $f_1(a)=e^{(a+1)x}=e^{ax}e^x.$ From
(\ref{f-e}), we derive $(qe^x+1)I_{-q}(e^{ax})=[2]_q.$ Hence we obtain
\begin{equation}\label{i-f-e}
I_{-q}(e^{ax})=\int_{\mathbb Z_p}e^{ax}d\mu_{-q}(a)=\frac{[2]_q}{qe^x+1}.
\end{equation}
We now set
\begin{equation}\label{q-e-d}
\frac{[2]_q}{qe^x+1}=\sum_{n=0}^\infty E^*_{n,q}\frac{x^n}{n!}.
\end{equation}
$E^*_{n,q}$ is called $q$-Euler numbers. By (\ref{i-f-e}) and (\ref{q-e-d}), we see that
$\int_{\mathbb Z_p}a^nd\mu_{-q}(a)=E^*_{n,q}.$
From (\ref{i-f-e}), we also note that
\begin{equation}\label{e-p}
\int_{\mathbb Z_p}e^{(t+a)x}d\mu_{-q}(a)=\frac{[2]_q}{qe^x+1}e^{tx}.
\end{equation}
In view of (\ref{q-e-d}) and (\ref{e-p}), we can consider $q$-Euler polynomials associated to $t$ as follows:
\begin{equation}\label{i-e-p}
\frac{[2]_q}{qe^x+1}e^{tx}=\sum_{n=0}^\infty E^*_{n,q}(t)\frac{x^n}{n!} \quad{\rm and}\quad \int_{\mathbb Z_p}(t+a)^nd\mu_{-q}(a)
=E^*_{n,q}(t).
\end{equation}
Put $\lim_{q\rightarrow1}E^*_{n,q}=E^*_n$ and $\lim_{q\rightarrow1}E^*_{n,q}(t)=E^*_n(t).$ Then we have $E_n(t)=E^*_n(t)$ and
$$E_n=\sum_{m=0}^n 2^m\binom nmE_m^*,$$
where $E_n$ and $E_{n}(t)$ are the ordinary
Euler numbers and polynomials.
By (\ref{q-e-d}) and (\ref{i-e-p}), we easily see that
$E^*_{n,q}(t)=\sum_{m=0}^n\binom nm t^{n-m}E^*_{m,q}.$
For $d\in\mathbb Z^+,$ let $f_d(a)=f(a+d).$ Then we have
\begin{equation}\label{q^nI}
q^dI_{-q}(f_d)+(-1)^{d-1}I_{-q}(f)=[2]_q\sum_{a=0}^{d-1}(-1)^{d-a-1}q^af(a),\quad
\text{see \cite{KKJR}.}
\end{equation}
If $d$ is odd positive integer, we
have
\begin{equation}\label{q^nI-o}
q^dI_{-q}(f_d)+I_{-q}(f)=[2]_q\sum_{a=0}^{d-1}(-1)^{a}q^af(a).
\end{equation}
Let $\chi$ be a Dirichlet character with conductor $d=d_\chi$(=odd)$\in\mathbb Z^+.$ If take $f(a)=\chi(a)e^{(t+a)x},$
then we have $f_d(a)=f(a+d)=\chi(a)e^{dx}e^{(t+a)x}.$ From (\ref{I_q}) and (\ref{q^nI-o}), we derive
\begin{equation}\label{i-ge-p}
\int_{X}\chi(a)e^{(t+a)x}d\mu_{-q}(a)=\frac{[2]_q\sum_{a=1}^{d}(-1)^aq^a\chi(a)e^{(t+a)x}}{q^de^{dx}+1}.
\end{equation}
In view of (\ref{i-ge-p}), we also consider the generalized $q$-Euler polynomials attached to $\chi$ as follows:
\begin{equation}\label{gen-p}
F_{\chi,q}(x,t)=\frac{[2]_q\sum_{a=1}^{d}(-1)^aq^a\chi(a)e^{(t+a)x}}{q^de^{dx}+1}=\sum_{n=0}^\infty E^*_{n,\chi,q}(t)\frac{x^n}{n!}.
\end{equation}
From (\ref{i-ge-p}) and (\ref{gen-p}), we derive the following
\begin{equation}\label{i-gen-p}
\int_{X}\chi(a)(t+a)^nd\mu_{-q}(a)=E^*_{n,\chi,q}(t)
\end{equation}
for $n\geq0.$ Put $\lim_{q\rightarrow1}E^*_{n,\chi,q}(t)=E^*_{n,\chi}(t).$
On the other hand, the generalized $q$-Euler polynomials attached to $\chi$ are easily expressed as the $q$-Euler polynomials:
\begin{equation}\label{r-gen-p}
E^*_{n,\chi,q}(t)=d^n\frac{[2]_q}{[2]_{q^d}}\sum_{a=1}^d(-1)^aq^a\chi(a)E^*_{n,q^d}\left(\frac{a+t}{d}\right),\quad n\geq0.
\end{equation}

Let $\chi$ be a Dirichlet character with conductor $d=d_\chi\in\mathbb Z^+.$ It is well known (see \cite{I, Wa1})
that, for positive integers $m$ and $n,$
\begin{equation}\label{p-s}
\sum_{a=1}^{dn}\chi(a)a^{m}=\frac1{m+1}(B_{m+1,\chi}(dn)-B_{m+1,\chi}(0)),
\end{equation}
where $B_{m+1,\chi}(t)$ is the generalized Bernoulli polynomials. When $d=d_\chi$(=odd)$\in\mathbb Z^+,$ note that
\begin{equation}\label{f-eq}
\begin{aligned}
&\frac{[2]_q\sum_{a=1}^{d}(-1)^aq^a\chi(a)e^{(t+a)x}(1-(-q^de^{dx})^n)}{1-(-q^de^{dx})} \\
&=[2]_q\sum_{a=1}^d\sum_{l=0}^{n-1}(-1)^{a+dl}q^{a+dl}\chi(a+dl)e^{x(t+a+dl)}\\
&=[2]_q\sum_{a=1}^{dn}(-1)^{a}q^{a}\chi(a)e^{x(t+a)} \\
&=\sum_{m=0}^\infty\left([2]_q\sum_{a=1}^{dn}(-1)^{a}q^{a}\chi(a)(t+a)^m\right)\frac{x^m}{m!}.
\end{aligned}
\end{equation}
By (\ref{gen-p}), the relation (\ref{f-eq}) can be rewritten as
\begin{equation}\label{f-e-g}
\begin{aligned}
&\frac{[2]_q\sum_{a=1}^{d}(-1)^aq^a\chi(a)e^{(t+a)x}(1-(-q^de^{dx})^n)}{1-(-q^de^{dx})} \\
&=\sum_{m=0}^\infty \left(E^*_{m,\chi,q}(t)+(-1)^{n+1}q^{dn}E^*_{m,\chi,q}(t+dn)\right)\frac{x^m}{m!}.
\end{aligned}
\end{equation}
Now, we give the $q$-analogue of (\ref{p-s}) for the generalized Euler polynomials. From (\ref{f-eq}) and (\ref{f-e-g}), it is easy to see that
\begin{equation}\label{e-s-p}
\sum_{a=1}^{dn}(-1)^{a}q^{a}\chi(a)(t+a)^m=\frac1{[2]_q}\left(E^*_{m,\chi,q}(t)+(-1)^{n+1}q^{dn}E^*_{m,\chi,q}(t+dn)\right)
\end{equation}
for positive integers $m$ and $n.$ In particular, replacing $q$ by 1 in (\ref{e-s-p}), if $\chi=\chi^0,$ the principal character $(d_\chi=1),$
and $t=0,$ then
$$\sum_{a=1}^{n-1}(-1)^a a^m=\frac12\left(E_m(0)+(-1)^{n+1}E_m(n)\right).$$

\begin{definition}\label{q-zeta}
Let $s\in\mathbb C$ with Re$(s)>1.$ Let $\chi$ be a primitive Dirichlet character with conductor $d=d_\chi$(=odd)$\in\mathbb Z^+.$
We set
$$l_q(s,t,\chi)=[2]_q\sum_{n=0}^\infty\frac{(-1)^nq^n\chi(n)}{(t+n)^s},\quad 0<t\leq1.$$
\end{definition}

\begin{remark}\label{q-zeta-in}
We assume that $q\in\mathbb C$ with $|q|<1.$
Let $\chi$ be a primitive Dirichlet character with conductor $d=d_\chi$(=odd)$\in\mathbb Z^+.$
From (\ref{gen-p}), we consider the below integral which known the Mellin transformation of
$F_{\chi,q}(x,t)$ (cf. \cite{Si}).
$$\begin{aligned}
\frac1{\Gamma(s)}\int_0^\infty x^{s-1}F_{\chi,q}(-x,t)dx
&=[2]_q \sum_{a=1}^{d}(-1)^aq^a\chi(a)\frac1{\Gamma(s)}\int_0^\infty x^{s-1}\frac{e^{-(t+a)x}}{1-(-q^de^{-dx})}dx \\
&=[2]_q \sum_{a=1}^{d}(-1)^aq^a\chi(a+dl)\sum_{l=0}^\infty(-1)^{dl}\frac{q^{dl}}{(a+dl+t)^{s}}.
\end{aligned}$$
We write $n=a+dl,$ where $n=1,2,\ldots,$ and obtain
$$\frac1{\Gamma(s)}\int_0^\infty x^{s-1}F_{\chi,q}(-x,t)dx=[2]_q\sum_{n=0}^\infty\frac{(-1)^nq^n\chi(n)}{(t+n)^s}=l_q(s,t,\chi).$$
\end{remark}

Note that $l_q(s,t,\chi)$ is an analytic function in the whole complex $s$-plane.

By using a geometric series in (\ref{gen-p}), we obtain
$$[2]_qe^{tx}\sum_{n=0}^\infty(-1)^nq^n\chi(n)e^{nx}=\sum_{n=0}^\infty E^*_{n,\chi,q}(t)\frac{x^n}{n!}.$$
We also note that
\begin{equation}\label{n-v}
E^*_{n,\chi,q}(t)=\left(\frac{\rm d}{{\rm d}x}\right)^k[2]_qe^{tx}\sum_{n=0}^\infty(-1)^nq^n\chi(n)e^{nx}\biggl|_{x=0}.
\end{equation}
By Definition \ref{q-zeta} and (\ref{n-v}), we obtain the following theorem.

\begin{proposition}
For $n\in\mathbb Z^+,$ we have $l_q(-n,t,\chi)=E^*_{n,\chi,q}(t).$
\end{proposition}

These values of $l_q(s,t,\chi)$ at netative integers are algebraic, hence may be regarded as being in an extension of $\mathbb Q_p.$
We therefore look for a $p$-adic function which agrees with $l_q(s,t,\chi)$ at the negative integers in Section \ref{p-adic-lq}.

\section{A two-variable $p$-adic $l_q$-function}\label{p-adic-lq}

We shall consider the $p$-adic analogue of the $l_q$-functions which are introduced in the previous section (see Definition \ref{q-zeta}).
Throughout this section we assume that $p$ is an odd prime. Note that there exists $\varphi(p)$ distinct solutions,
modulo $p,$ to the equation $x^{\varphi(p)}-1=0,$ and each solution must be congruent to one of the
values $a\in\mathbb Z,$ where $1\leq a< p,(a,p)=1.$ Thus, given $a\in\mathbb Z$ with $(a, p)=1,$ there exists a unique
$\omega(a)\in\mathbb Z_p,$ where $\omega(a)^{\varphi(p)}=1,$ such that $\omega(a)\equiv a\pmod{p\mathbb Z_p}.$
Letting $\omega(a)=0$ for $a\in\mathbb Z$ such that $(a,p)\neq1,$ it can be seen that $\omega$ is actually a Dirichlet character having conductor
$d_\omega=p,$ called the Teichm\"uller character. Let $\langle a\rangle=\omega^{-1}(a)a.$ Then $\langle a\rangle\equiv 1\pmod{p\mathbb Z_p}.$
For the context in the sequel, an extension of the definition of the Teichm\"uller character is
needed. We denote a particular subring of $\mathbb C_p$ as
$$R=\{ a\in\mathbb C_p\mid |a|_p\leq1\}.$$
If $t\in \mathbb C_p$ such that $|t|_p\leq1,$ then for any $a\in\mathbb Z, a+pt=a\pmod{p R}$ Thus, for $t\in\mathbb C_p,
|t|_p\leq1,\omega(a+pt)=\omega(a).$ Also, for these values of $t,$ let $\langle a+pt\rangle=\omega^{-1}(a)(a+pt).$
Let $\chi$ be the Dirichlet character of conductor $d=d_\chi.$ For $n\geq1,$ we define $\chi_n$ to be the primitive character
associated to the character $\chi_n:(\mathbb Z/\text{lcm}(d,p)\mathbb Z)^\times\rightarrow\mathbb C^\times$ defined by $\chi_n(a)=\chi(a)\omega^{-n}(a).$

\begin{definition}\label{def-lq}
Let $\chi$ be the Dirichlet character with conductor $d=d_\chi$(=odd) and let $F$ be a positive integral multiple
of $p$ and $d.$
Now, we define the two-variable $p$-adic $l_q$-functions as follows:
$$l_{p,q}(s,t,\chi)=\frac{[2]_q}{[2]_{q^F}}
\sum_{\substack{a=1\\ (p,a)=1}}^F(-1)^a\chi(a){q^a}\langle a+pt\rangle^{-s}\sum_{m=0}^\infty\binom {-s}m \left(\frac F{\langle a+pt\rangle}\right)^{m}
E^*_{m,q^F}.$$
\end{definition}
Let $D=\{s\in\mathbb C_p\mid |s|_p<p^{1-\frac1{p-1}}\}$ and
let $a\in\mathbb Z, (a,p)=1.$ For $t\in\mathbb C_p,|t|_p\leq1,$ the same argument as that given in the proof of the main theorem of
\cite{Fo, Wa2} can be
the functions $\sum_{m=0}^\infty\binom {s}m \left(F/(a+pt)\right)^{m}E^*_{m,q^F}$ and
$\langle a+pt\rangle^{s}=\sum_{m=0}^\infty\binom{s}{m}(\langle a+pt\rangle-1)^m$ is analytic for $s\in D.$
According to this method, we see that the function $\sum_{m=0}^\infty\binom {s}m \left(F/(a+pt)\right)^{m}E^*_{m,q^F}$ is analytic
for $t\in\mathbb C_p, |t|_p\leq1,$ whenever $s\in D.$ It readily follows that
$\langle a+pt\rangle^{s}=\langle a\rangle^{s}\sum_{m=0}^\infty\binom{s}{m}(a^{-1}pt)^m$ is analytic
for $t\in\mathbb C_p, |t|_p\leq1,$ when $s\in D.$
Therefore, $l_q(s,t,\chi)$ is analytic for $t\in\mathbb C_p, |t|_p\leq1,$ provided $s\in D$ (see \cite{Fo}).

We set
\begin{equation}\label{Hpq}
h_{p,q}(s,t,a|F)={(-1)^aq^a}\langle a+pt\rangle^{-s}
\frac{[2]_q}{[2]_{q^F}}\sum_{m=0}^\infty\binom {-s}m \left(\frac F{a+pt}\right)^{m}E^*_{m,q^F}.
\end{equation}
Thus,we note that
\begin{equation}\label{Hpq-n}
h_{p,q}(-n,t,a|F)=\omega^{-n}(a)(-1)^aq^aF^n\frac{[2]_q}{[2]_{q^F}}E^*_{n,q^F}\left(\frac{a+pt}{F}\right)
\end{equation}
for $n\in\mathbb Z^+.$ We also consider the two-variable $p$-adic $l_q$-functions which interpolate the generalized $q$-Euler polynomials at
negative integers as follows:
\begin{equation}\label{lpq}
l_{p,q}(s,t,\chi)=\sum_{\substack{ a=1\\ (p,a)=1}}^F\chi(a)h_{p,q}(s,t,a|F).
\end{equation}
We will in the process derive an explicit formula for this function. Before we begin this derivation, we need the following result concerning
generalized $q$-Euler polynomials:

\begin{lemma}\label{EF}
Let $F$ be a positive integral multiple of $d=d_\chi.$ Then for each $n\in\mathbb Z, n\geq0,$
$$E^*_{n,\chi,q}(t)=F^n\frac{[2]_q}{[2]_{q^F}}\sum_{a=1}^{F}(-1)^aq^a\chi(a)E^*_{n,q^F}\left(\frac{a+t}F\right).$$
\end{lemma}

We can derive by a manipulation of an appropriate generating functions.

Set $\chi_n=\chi\omega^{-n}.$ From (\ref{Hpq-n}) and (\ref{lpq}), we obtain
\begin{equation}\label{lpq-n}
\begin{aligned}
l_{p,q}(-n,t,\chi)
&=F^n\frac{[2]_q}{[2]_{q^F}}\sum_{\substack{ a=1\\(p,a)=1}}^F\chi_n(a)(-1)^aq^aE^*_{n,q^F}\left(\frac{a+pt}F\right)\\
&=F^n\frac{[2]_q}{[2]_{q^F}}\sum_{a=1}^F\chi_n(a)(-1)^aq^aE^*_{n,q^F}\left(\frac{a+pt}F\right) \\
&\qquad-F^n\frac{[2]_q}{[2]_{q^F}}\sum_{a=1}^{\frac Fp}\chi_n(pa)(-1)^{pa}q^{pa}E^*_{n,q^{F}}\left(\frac{pa+pt}F\right).
\end{aligned}
\end{equation}
for $n\in\mathbb Z^+.$ From Lemma \ref{EF}, we see that
\begin{equation}\label{e-int}
E^*_{n,\chi_n,q}(pt)=F^n\frac{[2]_q}{[2]_{q^F}}\sum_{a=1}^{F}(-1)^aq^a\chi_n(a)E^*_{n,q^F}\left(\frac{a+pt}F\right)
\end{equation}
and
\begin{equation}\label{e-pint}
E^*_{n,\chi_n,q^p}(t)=\left(\frac Fp\right)^n\frac{[2]_{q^p}}{[2]_{(q^{p})^{\frac Fp}}}
\sum_{a=1}^{\frac Fp}(-1)^a(q^{p})^a\chi_n(a)E^*_{n,(q^{p})^{\frac Fp}}
\left(\frac{a+t}{\frac Fp}\right).
\end{equation}

From (\ref{lpq-n}), (\ref{e-int}) and (\ref{e-pint}), we obtain the following theorem:

\begin{theorem}\label{lpq}
Let $F$(=odd) be a positive integral multiple of $p$ and $d=(d_\chi).$ Then
the two-variable $p$-adic $l_q$-functions
$$l_{p,q}(s,t,\chi)=\frac{[2]_q}{[2]_{q^F}}
\sum_{\substack{a=1\\ (p,a)=1}}^F(-1)^a\chi(a){q^a}\langle a+pt\rangle^{-s}\sum_{m=0}^\infty\binom {-s}m \left(\frac F{\langle a+pt\rangle}\right)^{m}
E^*_{m,q^F}$$
admits an analytic function for $t\in\mathbb C_p$ with $|t|_p\leq1$ and $s\in D,$ and satisfies the relation
$$l_{p,q}(-n,t,\chi)=E^*_{n,\chi_n,q}(pt)-p^n\chi_n(p)\frac{[2]_q}{[2]_{q^p}}E^*_{n,\chi_n,q^p}(t)$$
for $n\in\mathbb Z^+$ and $t\in\mathbb C_p$ with $|t|_p\leq1$
\end{theorem}

From (\ref{lpq}) and Theorem \ref{lpq}, it follows that $h_{p,q}(s,t,a|F)$ is analytic for $t\in\mathbb C_p$ with $|t|_p\leq1$ and $s\in D.$

\begin{remark} Let $\langle a+pt\rangle=\omega^{-1}(a)(a+pt),$ and let $t\in\mathbb C_p$ with $|t|_p\leq1$ and $s\in D.$ Then
the two-variable $p$-adic $l_q$-functions defined above is redefined by
$$l_{p,q}(s,t,\chi)=\int_{X^*}\chi(a)\langle a+pt\rangle^{-s} d\mu_{-q}(a),\quad \text{cf. \cite{KKJR, Yo}}.$$
Then we have
$$\begin{aligned}
l_{p,q}(-n,t,\chi)&=\int_{X}\chi_n(a)(a+pt)^n d\mu_{-q}(a)-\int_{X}\chi_n(pa)(pa+pt)^n d\mu_{-q}(pa) \\
&\overset{(\ref{i-gen-p})}{=}E^*_{n,\chi_n,q}(pt)-p^n\chi_n(p)\frac{[2]_q}{[2]_{q^p}}E^*_{n,\chi_n,q^p}(t),
\end{aligned}$$
since $X^*=X-pX$ and $[2]_{q^p}d\mu_{-q}(pa)=[2]_qd\mu_{-q^p}(a).$
\end{remark}

\begin{corollary} Let $F$(=odd) be a positive integral multiple
of $p$ and $d=(d_\chi),$ and let the two-variable $p$-adic $l$-functions
$$l_{p}(s,t,\chi)=
\sum_{\substack{a=1\\ (p,a)=1}}^F(-1)^a\chi(a)\langle a+pt\rangle^{-s}\sum_{m=0}^\infty\binom {-s}m \left(\frac F{\langle a+pt\rangle}\right)^{m}
E^*_{m}.$$
Then
\begin{enumerate}
\item $l_{p}(s,t,\chi)$ is analytic for $t\in\mathbb C_p$ with $|t|_p\leq1$ and $s\in D.$
\item $l_{p}(-n,t,\chi)=E^*_{n,\chi_n}(pt)-p^n\chi_n(p)E^*_{n,\chi_n}(t)$ for $n\in\mathbb Z^+.$
\item $l_{p}(s,t,\chi)=\int_{X^*}\chi(a)\langle a+pt\rangle^{-s} d\mu_{-1}(a)$ for $t\in\mathbb C_p$ with $|t|_p\leq1$ and $s\in D.$
\end{enumerate}
\end{corollary}

\end{document}